\documentclass[1p]{elsarticle}
\usepackage{amssymb}
\usepackage{amsfonts}

\newtheorem{theorem}{Theorem}
\newtheorem{conjecture}[theorem]{Conjecture}

\newtheorem{lemma}[theorem]{Lemma}
\newtheorem{observation}[theorem]{Observation}

\newproof{pf}{Proof}

\begin{document}
\title{On the Standard $(2,2)$-Conjecture}

\author[agh]{Jakub Przyby{\l}o\fnref{MNiSW}}

\fntext[MNiSW]{This work was partially supported by the Faculty of Applied Mathematics AGH UST statutory tasks within subsidy of Ministry of Science and Higher Education.}

\address[agh]{AGH University of Science and Technology, Faculty of Applied Mathematics, al. A. Mickiewicza 30, 30-059 Krakow, Poland}

\begin{abstract}
The well-known 1--2--3 Conjecture asserts that the edges of every graph without an isolated edge can be weighted with $1$, $2$ and $3$ so that adjacent vertices receive distinct weighted degrees. This is open in general. We prove that every graph with minimum degree $\delta\geq 10^6$ can be decomposed into two subgraphs requiring just weights $1$ and $2$ for the same goal. We thus prove the so-called Standard $(2,2)$-Conjecture for graphs with sufficiently large minimum degree. The result is in particular based on applications of the Lov\'asz Local Lemma and theorems on degree-constrained subgraphs.
\end{abstract}


\maketitle

\section{Introduction}

In 2004 Karo\'nski, {\L}uczak and Thomason~\cite{123KLT} posed a basic and apparently inconspicuous question. 
It is nowadays usually referred to as the \emph{1--2--3 Conjecture}. 
\begin{conjecture} 
Every connected graph with at least three vertices can be edge weighted with $1,2,3$ so that adjacent vertices receive distinct sums of their incident weights. 
\end{conjecture}
This occurred to be highly non-trivial.
The authors of~\cite{123KLT} were not in particular able to settle a constant upper bound concerning this concept, that is to verify if there is a constant $K$ such that each  connected graph $G=(V,E)$ with $|V|\geq 3$ is $\{1,2,\ldots,K\}$-\emph{weight colourable}, i.e. that there is a function $c:E\to\{1,2,\ldots,K\}$ such that the \emph{weighted degree} of a vertex $v$:
$$d_c(v):=\sum_{u\in N(v)}c(uv)$$ 
is distinct from the weighted degrees of all its neighbours for every $v\in V$ -- such $c$ is 
called a \emph{neighbour sum distinguishing $k$-weighting} (or \emph{colouring}) of $G$. 
They however showed that such a weighting always exists from a finite set $C$ of real numbers -- this was further investigated by Addario-Berry et al.~\cite{Louigi2}, and quite recently by Vu\v{c}kovi\'c~\cite{Vuckovic_3-multisets}, who finally settled its minimal size at $3$.  
The 1--2--3 Conjecture itself has gained considerable interest within combinatorial community, see e.g.~\cite{Louigi30,Louigi,Julien5regular123,DudekWajc123complexity,KalKarPf_123,1234Reg123,Seamon123survey,ThoWuZha,123with13}
 (and~\cite{BarGrNiw,PrzybyloWozniakChoos,WongZhu23Choos,WongZhuChoos} for results concerning related concepts in a list setting), yet still remains open.
In~\cite{Louigi30} Addario-Berry et al. showed the first constant upper bound supporting it by proving that a set of 
$30$ least positive integers is sufficient. The size of this was then pushed down to $16$ by Addario-Berry, Dalal and Reed~\cite{Louigi}, and next to $13$ by Wang and Yu~\cite{123with13}. The best general result thus far asserts that neighbours in every connected graph with at least three vertices can be sum-distinguished by means of weights $1,2,3,4,5$, see~\cite{KalKarPf_123} by Kalkowski, Karo\'nski and Pfender.
Moreover, every $d$-regular graph, $d\geq 2$ is known to be $\{1,2,3,4\}$-weight colurable,~\cite{1234Reg123}, while asymptotically almost surely a random graph ($G_{n,p}$ with constant
$p\in(0,1)$) is $\{1,2\}$-weight colurable,~\cite{Louigi}, even though
Dudek and Wajc~\cite{DudekWajc123complexity} showed that determining whether a particular graph 
is $\{1,2\}$-weight colurable is NP-complete.  
Many of these results were obtained due to development and application of a very convenient tool concerning sufficient conditions for existence subgraphs with some desired 
 features in a given graph, so-called \emph{degree constrained subgraphs} (see e.g.~\cite{Louigi30,Louigi2,Louigi}). We shall also make use of such a handy result, see Lemma~\ref{1_6Lemma} below.

The 1--2--3 Conjecture can alternatively be expressed in terms of \emph{lcocally irregullar multigraphs}, which we understand here as multigraphs with distinct degrees of adjacent vertices. In such a setting we simply ask if the edges of a connected graph with at least three vertices can be multiplied, each at most $3$ times (counting in the original copy of an edge) so that the obtained multigraph is locally irregular.
An edge colouring of a given graph can be equivalently viewed at as its (edge) decomposition into subgraphs -- each induced by edges of a single colour.
In~\cite{LocalIrreg_1} Baudon et al. started research on decomposability of graphs into locally irregular sub\emph{graphs}, and posed a conjecture that except for some family $\mathfrak{T}$ of specific graphs (of maximum degree at most three), all connected graphs can be decomposed into $3$ locally irregular subgraphs. 
In fact this problem is strongly motivated by the 1--2--3 Conjecture and widely interrelated with this concept (note in particular that a locally irregular graph fulfills the 1--2--3 Conjecture, it is even $\{1\}$-weight colourable), 
see~\cite{LocalIrreg_1} for details.  In~\cite{LocalIrreg_2} the conjecture from~\cite{LocalIrreg_1} was confirmed for graphs with minimum degree at least $10^{10}$, while a constant upper bound in case of all graphs (except those in $\mathfrak{T}$) was provided by Bensmail, Merker and Thomassen~\cite{BensmailMerkerThomassen}, and then slightly optimized by Lu\v{z}ar, Przyby{\l}o, and Sot\'ak~\cite{LocalIrreg_Cubic}.

In this paper we prove that any graph $G$ can be decomposed into two $\{1,2\}$-weight colourable subgraphs if only the minimum degree of $G$ is large enough, see Theorem~\ref{MainTh22Standard} below. (Note that the mentioned random graph $G_{n,p}$ with constant $p\in(0,1)$, which is asymptotically almost surely $\{1,2\}$-weight colourable is very likely ``to be close'' to a regular graph, and thus often admits uncomplicated and straightforward application of probabilistic tools such as the Lov\'asz Local Lemma -- the lack of such a convenience was one of the main obstacles we had to overcome within our random approach.)
Thereby we significantly improve the result of~\cite{BensmailPrzybylo_Decomposability23}, which implies in particular that every graph $G$ without isolated edges can be decomposed into 24 $\{1,2,3\}$-weight colourable subgraphs (i.e. subgraphs fulfilling the 1--2--3 Conjecture), or at most $2$ such subgraphs if $G$ is $d$-regular, $d\geq 18$.
We thus also show that the following so-called \emph{Standard $(2,2)$-Conjecture} from~\cite{8Authors} holds for graphs with sufficiently large minimum degree.
\begin{conjecture}\label{Standard22Conjecture-quoted}
Every graph without isolated edges and isolated triangles
can be decomposed into two $\{1,2\}$-weight colourable subgraphs.
\end{conjecture}
Research concerning this conjecture were naturally motivated by the concepts discussed above.
This paper also improves the result associated with the weaker version of Conjecture~\ref{Standard22Conjecture-quoted}, the so-called  \emph{Weak $(2,2)$-Conjecture}, which postulates that there exists a decomposition of any graph $G=(V,E)$ without isolated edges into two subgraphs $G_1,G_2$ and corresponding $\{1,2\}$-weightings $c_1,c_2$ of these so that $d_{c_1}(u)\neq d_{c_1}(v)$ or $d_{c_2}(u)\neq d_{c_2}(v)$
for every $uv\in E$ (note this condition is implied by the one within the Standard $(2,2)$-Conjecture -- see~\cite{8Authors} and~\cite{Weak22Conjecture} for details), which in~\cite{Weak22Conjecture} was proved to hold for graphs with large enough minimum degree.

\section{Basic Tools}

We first present one basic observation followed by a recollection of a few fundamental tools of 
the probabilistic method we shall use later on.  
\begin{observation}
\label{EvenDecompositionOptimal}
Every connected graph $G=(V,E)$ with minimum degree $\delta$ can be decomposed into two subgraphs $G_1, G_2$ so that for every $v\in V$ and $i\in{1,2}$,
\begin{equation}\label{EQ_EulerianDecomposition}
d_{G_i}(v)\geq \left\lfloor \frac{d_G(v)}{2}\right\rfloor,
\end{equation} 
except possibly one vertex $u\in V$ with $d_G(u)=\delta$ for which we may have 
$d_{G_1}(u)=\lceil\frac{\delta+1}{2}\rceil$ (and hence $d_{G_2}(u)=\lfloor\frac{\delta-1}{2}\rfloor$).
\end{observation}
\begin{pf}
If $G=(V,E)$ contains vertices of odd degree, we join them by edges with a single new vertex $u$ and denote the obtained graph by $H$; otherwise we set $H:=G$ and denote any vertex with minimum degree in $G$ as $u$. Obviously all degrees of $H$ are even then, and thus there exists an Eulerian tour in it.   We then fix a decomposition of $H$ into $H_1$ and $H_2$ by traversing the Eulerian tour starting from $u$ and alternately including the encountered edges in $H_1$ and $H_2$. Finally we define $G_1=H_1$ and $G_2=H_2$ if $G$ contained no odd degree vertices, and set $G_1=H_1-u$, $G_2=H_2-u$ otherwise. It is straightforward to verify that in the latter case the inequality (\ref{EQ_EulerianDecomposition}) holds for all vertices $v\in V$ ($i=1,2$), while otherwise (\ref{EQ_EulerianDecomposition}) may not hold for only one vertex, namely $u$, for which we then have $d_{G_1}(u)=\lceil\frac{\delta+1}{2}\rceil$. 
\qed
\end{pf}

The following standard version of the Lov\'asz Local Lemma can be found e.g. in~\cite{MolloyReed}. 

\begin{theorem}[\textbf{The Local Lemma}]
\label{LLL-general}
Let $\mathcal{A}$ be a finite family of events in any probability space and let $D=(\mathcal{A},E)$ be a directed graph such that every event $A\in \mathcal{A}$ is mutually independent of all the events $\{B\in\mathcal{A}: (A,B)\notin E, B\neq A\}$.
Suppose that there are real numbers $p_A$ ($A\in\mathcal{A}$) such that for every $A\in\mathcal{A}$, $0\leq p_A<1$ and
\begin{equation}\label{EqLLL-general}
{\rm \emph{\textbf{Pr}}}(A) \leq p_A \prod_{B\leftarrow A} (1-p_B).
\end{equation}
Then ${\rm \emph{\textbf{Pr}}}(\bigcap_{A\in\mathcal{A}}\overline{A})>0$.
\end{theorem}
Here $B\leftarrow A$ (or $A\rightarrow B$) means that there is an arc from $A$ to $B$ in $D$,
the so-called \emph{dependency digraph}.
The Chernoff Bound below can be found e.g. in~\cite{JansonLuczakRucinski} (Th. 2.1, page 26).

\begin{theorem}[\textbf{Chernoff Bound}]\label{ChernofBoundTh}
For any $0\leq t\leq np$,
$${\mathbf Pr}({\rm BIN}(n,p)>np+t)<e^{-\frac{t^2}{3np}}~~~~{and}~~~~{\mathbf Pr}({\rm BIN}(n,p)<np-t)<e^{-\frac{t^2}{2np}}$$ 
where ${\rm BIN}(n,p)$ is the sum of $n$ independent Bernoulli variables, each equal to $1$ with probability $p$ and $0$ otherwise.
\end{theorem}

One more crucial in our argumentation lemma, concerning degree-constrained subgraphs is included in the next section, followed by a discussion on its utility in our setting.

\section{General Idea of Proof}

Our main result is the following, cf. Conjecture~\ref{Standard22Conjecture-quoted}.
\begin{theorem}\label{MainTh22Standard}
Every graph with minimum degree $\delta\geq 10^6$ 
can be decomposed into two $\{1,2\}$-weight colourable subgraphs.
\end{theorem}

In order to prove it we shall make use of the following lemma, which can be found in~\cite{LocalIrreg_2}, and which is a direct corollary of a theorem from~\cite{Louigi30}
(see also \cite{Louigi2,Louigi} for similar theorems on degree-constrained subgraphs and their applications).

\begin{lemma} \label{1_6Lemma}
Suppose that for some graph $G=(V,E)$ with minimum degree at least $12$ we have chosen, for every
vertex $v$, an integer $\lambda_v\geq 2$ with $6\lambda_v \leq d(v)$. Then for every
assignment
$$a:V\to \mathbb{Z},$$
there exists a spanning subgraph $H$ of $G$ such that $d_H(v)\in[\frac{d(v)}{3},\frac{2d(v)}{3}]$ and
$d_H(v)\equiv a(v) \pmod {\lambda_v}$ or $d_H(v)\equiv a(v)+1 \pmod {\lambda_v}$ for each $v\in V$.
\end{lemma}

In order to exemplify its practical utility, let us consider a graph $G=(V,E)$ with minimum degree $\delta\geq 12$, and suppose it has a relatively small chromatic number, that is $\chi(G)\leq \delta/12$. Let $a':V\to\{0,2,4,\ldots, 2\chi(G)-2\}$ be a proper vertex colouring of $G$, and let us apply Lemma~\ref{1_6Lemma} to $G$ with an assignment $a:V\to \mathbb{Z}$ such that $a(v)=a'(v)-d(v)$ and with $\lambda_v=2\chi(G)$ for every $v\in V$. Then we define a $\{1,2\}$-weighting $c$ by attributing weight $2$ to all edges of the resulting subgraph $H$, and weight $1$ to the remaining edges of $G$. Consequently, for every $v\in V$, $d_c(v)=2d_H(v)+1(d(v)-d_H(v))=d(v)+d_H(v) \equiv a'(v),a'(v)+1 \pmod {2\chi(G)}$. By our construction we thus have that $d_c(u)\neq d_c(v)$ for every $uv\in E$.

The  argument above thus implies that a graph is $\{1,2\}$-weight colourable provided that its minimum degree $\delta\geq 12$ is sufficiently many times (at least $12$) larger than its chromatic number. (Such conclusion was in fact already applied in~\cite{Louigi} for the case of random graphs.)
It would thus be convenient towards proving Theorem~\ref{MainTh22Standard} if we could provide a decomposition
of a given graph $G$ into two subgraphs with the chromatic number much smaller than the minimum degree.
Unfortunately we are not  in general able to achieve such a goal, even for graphs with large minimum degree. We thus shall construct a decomposition complying with weaker requirements, but still sufficient for our purposes. 
The foundations of our proof of Theorem~\ref{MainTh22Standard} are however based on the idea presented above, and elaborate the concept of accomplishing a desired decomposition in the case of regular graphs (with degree large enough) we outline below in this section.
In fact the primary idea behind such a decomposition of a regular graph in turn stems  from analysis of complete graphs.

For simplicity consider first a complete graph $K_{n^2}$ where $n$ is a (large enough) \emph{even} integer.
Let $V(K_{n^2})=\{0,1,2,\ldots,n-1\}\times \{0,1,2,\ldots,n-1\}$ (or in other words let us colour or label the vertices of 
$K_{n^2}$ with distinct ordered pairs of nonnegative integers not exceeding $n-1$). Then we define a decomposition of $K_{n^2}$ into two subgraphs $G_1=(V,E_1)$ and $G_2=(V,E_2)$ by including in $E_1$ all edges $(i,j)(k,l)$ such that $j=l$ and all with $i+j+k+l\equiv 1~({\rm mod}~2)$ and $i\neq k$. The remaining edges $(i,j)(k,l)$ of $K_{n^2}$ are included in $E_2$ -- note that either $i=k$ or  $i+j+k+l\equiv 2~({\rm mod}~2)$ and $j\neq l$ for these. 
It is not hard to verify that then:
$$\delta(G_1),\delta(G_2)\geq \left\lfloor\frac{n^2-1}{2}\right\rfloor\geq \frac{n^2}{2}-1.$$ 
On the other hand, by the definition of $G_1$, we may properly colour its vertices by choosing colour $i$ for every vertex $(i,j)\in V$ (vertices with the same first coordinates can be coloured the same, as they form an independent set in $G_1$). Analogously, we may colour properly $G_2$ by choosing colour $j$ for every vertex $(i,j)\in V$, and hence: 
$$\chi(G_1),\chi(G_2)\leq n.$$
Therefore, already for relatively small $n$, we shall have $\chi(G_i)\leq\delta(G_i)/12$ for $i=1,2$, and thus by the argument presented above (based on application of Lemma~\ref{1_6Lemma}), each of $G_1$ and $G_2$ 
shall be $\{1,2\}$-weight colourable.
(A similar reasoning can also be applied when the order of a complete graph is not a square of an even integer).

Now suppose we consider a $d$-regular graph $G=(V,E)$ with (large enough) $d$ instead of $K_{n^2}$, and we randomly assign colours in $\{0,1,2,\ldots,n-1\}\times \{0,1,2,\ldots,n-1\}$ (for some even $n\geq \sqrt{d}$) to the vertices of $G$. If the obtained colouring is proper (we describe later how to actually go around this particular issue) then it defines a preserving colours homomorphism of $G$ into $K_{n^2}$ coloured as described in the paragraph above, and we may define a decomposition of $G$ into $G_1=(V,E_1)$ and $G_2=(V,E_2)$ 
following the same rules as in the case of $K_{n^2}$. As the random distribution of colours ought to be relatively uniform, with high probability, the share of edges of $E_1$ and $E_2$ around every vertex -- mirroring the one in $K_{n^2}$ should be fairly proportional. Therefore the minimum degrees of $G_1$ and $G_2$ should not be much smaller than $d/2$, while by the construction the chromatic number of each of these shall be at most $n$ (analogously as in the case of a decomposition of $K_{n^2}$). Then if only $n$ is significantly smaller than $d/2$ (we in general certainly need roughly at least $n\geq \sqrt{d}$, in order to facilitate some approximation of properness of the vertex colouring of $G$ resulting from the random process), with positive probability  we should  
obtain a desired decomposition of $G$, whose validity shall follow by Lemma~\ref{1_6Lemma} as above.
We omit details concerning regular graphs here, as we intend to focus on the more demanding general case. 
Aiming at this we shall have to face several problems,
as e.g. the mentioned earlier fact that the Local Lemma is often unwieldy without any assumption on the relation between the maximum and minimum degree in a graph at hand.
In particular, while randomly assigning colours and weighted degrees in the following proof 
of Theorem~\ref{MainTh22Standard}, we shall need several refinements,  among others: 
\begin{itemize}
\item vertices with smaller degrees shall be granted (in some sense) a limited list of choices for their potential weighted degrees in $G_1$ and $G_2$, proportional to their degree (in $G$); these shall be proportionally augmented in the case of vertices of larger degrees;
\item in order to optimize the final result, the number $[n(v)]^2$ of admitted choices of colours for every vertex $v$ shall be of order $[d(v)]^2$ (rather than $[\sqrt{d(v)}]^2$), more precisely   a constant times smaller than $[d(v)]^2$;
\item by our construction, only the neighbours of relatively close degrees in $G$ (one at most twice as large as the other) shall present a potential threat of conflict between weighted degrees in $G_1$ or $G_2$, what shall facilitate application of the Lov\'asz Local Lemma;
\item the auxiliary vertex colouring we provide within the proof shall not be required to be proper -- we extend the lists of potential weighted degrees to overcome this problem instead 
(this is actually optional, but yields a better final result).
\end{itemize}

\section{Proof of Theorem~\ref{MainTh22Standard}}

\subsection{Basic Quantities}
 
Let $G=(V,E)$ be a graph with minimum degree $\delta\geq 10^6$.
In what follows, by $d(v)$ we shall always mean the degree of a vertex $v$ in $G$, i.e. $d_G(v)$.

Let $H'$ be the subgraph of $G$ induced by all the edges $uv\in E$ such that $d(u)\notin [d(v)/2,2d(v)]$.
(By our construction there shall be no possible sum-conflict between such $u$ and $v$ in the two 
$\{1,2\}$-weighted subgraphs of $G$ we are aiming to construct.) 
By Observation~\ref{EvenDecompositionOptimal} we may decompose $H'$ into two subgraphs $H'_1$ and $H'_2$ including almost equal share of incident edges of every vertex in $H'$, i.e. such that for each $v\in V$,
\begin{equation}\label{dH'iBound}
d_{H'_1}(v),d_{H'_2}(v)\geq \left\lfloor\frac{d_{H'}(v)-1}{2}\right\rfloor \geq \frac{d_{H'}(v)}{2}-1.
\end{equation}

Let then $H=G-E(H')$, i.e. 
\begin{equation}\label{H-degrees}
d(u)\in \left[\frac{d(v)}{2},2d(v)\right] 
~~~~{\rm for~ every}~~~~ uv\in E(H),
\end{equation}
and set
$$q:=\frac{9}{20}=0.45~~~~~~{\rm and}~~~~~~t:=18.$$
We shall randomly decompose $H$ to $H_1,H_2$ such that the two subgraphs of the obtained decomposition of $G$:
$$G_1:=H'_1\cup H_1,~~~~G_2:=H'_2\cup H_2$$
divide fairly equally incident edges of all the vertices, i.e. such that for every $v\in V$:
\begin{equation}\label{evenG1G2}
d_{G_1}(v),d_{G_2}(v)\geq qd(v). 
\end{equation}
Within the same random process we shall at the same time provide desired $\{1,2\}$-weightings of $G_1,G_2$.
For this aim we first define for each $v\in V$ a special even integer 
(corresponding in some sense to $n$ in the example concerning $K_{n^2}$ above):
$$
y_v:=2^{\lfloor\log_2\frac{qd(v)}{24 t}\rfloor}.
$$
Note that 
\begin{equation}\label{yvBound}
\frac{qd(v)}{24 t}\geq y_v\geq\frac{qd(v)}{48 t}
\end{equation}
and
\begin{equation}\label{yuyvBound}
y_u\in \left\{\frac{y_v}{2},y_v,2y_v\right\}~~~~{\rm for~every}~~~~uv\in E(H).
\end{equation}

\subsection{Random Assignment}

We first randomly and independently assign to every vertex $v\in V$
a pair of integers $(c^1_v,c^2_v)\in [0,y_v-1]^2$ -- each with equal probability.

Let $T$ be (a random variable expressing) 
the least integer such that every vertex $v\in V$ has at most $2T-2$ neighbours $u$ in $H$ with $y_u=y_v$ and $(c^1_u,c^2_u)=(c^1_v,c^2_v)$ 
(such problematic neighbours, coloured the same, shall later require blowing up lists of admissible weighted degrees for $v$; 
we shall distribute edges joining $v$ with them fairly evenly between $H_1$ and $H_2$ though).

We  define the edge sets of $H_1$ and $H_2$ as follows. Let $uv\in E(H)$. 
\begin{itemize}
\item[($1^\circ$)] if $c^1_u=c^1_v$ and $c^{2}_u\neq c^{2}_v$, then $uv\in E(H_2)$;
\item[($2^\circ$)] if $c^2_u=c^2_v$ and $c^{1}_u\neq c^{1}_v$, then $uv\in E(H_1)$;
\item[($3^\circ$)] if $c^1_u\neq c^1_v$, $c^{2}_u\neq c^{2}_v$ and $c^1_u+c^{2}_u+c^1_v+c^{2}_v \equiv 1~({\rm mod}~2)$, then $uv\in E(H_1)$;
\item[($4^\circ$)] if $c^1_u\neq c^1_v$, $c^{2}_u\neq c^{2}_v$ and $c^1_u+c^{2}_u+c^1_v+c^{2}_v \equiv 2~({\rm mod}~2)$, then $uv\in E(H_2)$;
\item[($5^\circ$)] for all non-negative integers $c^1,c^2,y$, let $H_{c^1,c^2,y}$ be the subgraph of $H$ induced by all vertices $v$ with $y_v=y$ and $(c^1_v,c^2_v)=(c^1,c^2)$;  we apply Observation~\ref{EvenDecompositionOptimal} to decompose it (one component after another) to $H^1_{c^1,c^2,y}$ and $H^2_{c^1,c^2,y}$ with all vertices of degree at most $T-1$ except at most one \emph{special} vertex in each component, which might be of degree $T$ -- we denote the set of all such special vertices by $V^*$; 
we include the edges of all $H^1_{c^1,c^2,y}$ in $E(H_1)$ (in case of existence of more than one such decomposition, we choose any of these via arbitrary deterministic rule fixed prior lunching the random process);
\item[($6^\circ$)]  we include in $E(H_2)$ the set of all the remaining edges of $H$ (comprising edges of all  $H^2_{c^1,c^2,y}$ and other edges $uv\in E(H)$ with $c^1_u = c^1_v$ and $c^{2}_u = c^{2}_v$, i.e. not belonging to any $H^1_{c^1,c^2,y}$ or $H^2_{c^1,c^2,y}$, hence with $y_u\neq y_v$).
\end{itemize}

Let us now set:
$$E_0:=\left\{uv\in E(H): (c^1_u,c^2_u) = (c^1_v,c^2_v)\right\}.$$
In order to avoid some potential dependences influencing application of the Local Lemma, instead of focusing on $H_1$ and $H_2$, we shall consider degrees in the following subgraphs of these:
$$H''_1:=H_1-E_0~~~~~~{\rm and}~~~~~~H''_2:=H_2-E_0,$$
for which the following obviously holds for every $v\in V$:
\begin{equation}\label{dHadH''}
d_{H_1}(v)\geq d_{H''_1}(v)~~~~~~{\rm and}~~~~~~d_{H_2}(v)\geq d_{H''_2}(v).
\end{equation}

Note that by ($1^\circ$) -- ($4^\circ$), similarly as in the case of $K_{n^2}$ above, regardless of the choice of the pair $(c^1_v,c^2_v)$ for a vertex $v\in V$, for every its neighbour $u$ in $H$ there are at least 
$$\frac{y_u^2-y_u}{2}$$ 
pairs out of all possible $y_u^2$ that might be assigned to $u$ which result in $uv\in E(H''_1)$
(note that unlike in the case of complete or more generally regular graphs we may have $c^1_v\notin[0,y_u-1]$ or $c^2_v\notin[0,y_u-1]$ this time, when $y_u<y_v$)
and analogously at least the same number of pairs causing $uv\in E(H''_2)$.
Therefore, by~(\ref{yvBound}), the probability that any given such $u$ (regardless of the choice for $v$) is assigned a pair leading to $uv\in E(H''_1)$ (or analogously $E(H''_2)$) equals at least (cf.~(\ref{yvBound})):
\begin{equation}\label{ProbabilEdge12Bound}
\frac{\frac{y_u^2-y_u}{2}}{y_u^2} =\frac12-\frac{1}{2y_u} \geq \frac12-\frac{48 t}{2qd(u)} \geq 
\frac12-\frac{24 \cdot 40}{10^6} \geq 0.49.
\end{equation}

We aim at showing that (with high enough probability) the edges of $G$ can be fairly evenly distributed between $G_1$ and $G_2$, i.e. that~(\ref{evenG1G2}) holds. For this goal we shall guarantee that
\begin{equation}\label{dHiBound}
d_{H''_i}(v) \geq qd_H(v)+1~~~~{\rm for}~~~~i=1,2
\end{equation}
for every vertex $v\in V$ with $d_H(v) > (1-2q)d(v)-2$, as then by~(\ref{dH'iBound}) and~(\ref{dHadH''}) we shall have
\begin{eqnarray}
d_{G_i}(v) &=& d_{H'_i}(v)+d_{H_i}(v) \geq d_{H'_i}(v)+d_{H''_i}(v) \nonumber\\
&\geq& \left(\frac{1}{2}d_{H'}(v)-1\right) + \left(qd_H(v)+1\right)\geq q\left(d_{H'}(v)+d_{H}(v)\right)=qd(v)\nonumber
\end{eqnarray}
for $i=1,2$. Note on the other hand that if $d_H(v) \leq (1-2q)d(v)-2$, then by~(\ref{dH'iBound}), 
$$d_{G_i}(v)\geq d_{H'_i}(v) \geq \frac{d_{H'}(v)}{2}-1 = \frac{d(v)-d_{H}(v)}{2}-1 \geq \frac{2qd(v)+2}{2}-1 = qd(v)$$
for $i=1,2$. In both cases~(\ref{evenG1G2}) shall hold.

The second feature we shall require from the sought assignments of pairs to the vertices is that
the problematic edges, those in all $H_{c^1,c^2,y}$, are not so frequent around any vertex, i.e. that $T\leq t$.
We thus show in the next subsection that the probability that some of our requirements does not hold for a given vertex is relatively small.

\subsection{Bad Events}

For any given vertex $v\in V$, let 
$$A(v)=\{u\in N_H(v): y_u=y_v \wedge (c_u^1,c_u^2)=(c_v^1,c_v^2)\},$$
and denote the following event
$$A_v:~~|A(v)|\geq 2t-1.$$
(Note that $T\leq t$ if and only if $A_v$ does not hold for any $v\in V$.)
As for any $u\in N_H(v)$ with $y_u=y_v$, the probability that $(c_u^1,c_u^2)=(c_v^1,c_v^2)$ equals exactly 
$y_u^{-2}=y_v^{-2}$ and the choices for  all vertices are independent, by~(\ref{yvBound}) we have: 
\begin{eqnarray}
\mathbf{Pr}(A_v) &\leq& {d(v) \choose 2t-1} \left(\frac{1}{y_v^2}\right)^{2t-1} \nonumber\\
&\leq& \frac{(d(v))^{2t-1}}{(2t-1)!} \cdot \frac{(48t)^{4t-2}}{(qd(v))^{4t-2}}\nonumber\\ 
&=&  \frac{1920^{70}}{(35!)(d(v))^{35}}. \label{Avineq} 
\end{eqnarray}

Next, for every $v\in V$, we denote the event:
$$B_v:~~d_H(v) > (1-2q)d(v)-2 ~\wedge~ \left(d_{H''_1}(v)< qd_H(v)+1 ~\vee~ d_{H''_2}(v)< qd_H(v)+1\right).$$
(Note that (\ref{dHiBound}) holds for every vertex $v\in V$ with $d_H(v) > (1-2q)d(v)-2$    if and only if   $B_v$ does not hold for all $v\in V$.) 
If $d_H(v) \leq (1-2q)d(v)-2$, then obviously
$$\textbf{Pr}\left(B_v\right)=0.$$
Suppose thus that $d_H(v) > (1-2q)d(v)-2=0.1d(v)-2$. Then by~(\ref{ProbabilEdge12Bound}) and the Chernoff Bound, for every pair of integers $(c^1,c^2)\in [0,y_v-1]^2$:
\begin{eqnarray}
&&\mathbf{Pr}\left(B_v ~|~ (c^1_v,c^2_v) = (c^1,c^2)\right)\nonumber\\  
&\leq& \mathbf{Pr}\left(d_{H''_1}(v)< qd_H(v)+1 ~\vee~ d_{H''_2}(v)< qd_H(v)+1 ~|~ (c^1_v,c^2_v) = (c^1,c^2)\right)\nonumber\\
&\leq& 2\cdot \mathbf{Pr}\left({\rm BIN}\left(d_H(v),0.49\right)< 0.45d_H(v)+1\right)\nonumber\\
&<& 2e^{-\frac{(0.04d_H(v)-1)^2}{2\cdot0.49d_H(v)}}\nonumber\\
&\leq& 2e^{-\frac{(0.04(0.1d(v)-2)-1)^2}{2\cdot0.49(0.1d(v)-2)}}\nonumber\\
&\leq& 2e^{-\frac{(0.0035d(v))^2}{0.1d(v)}}\nonumber\\
&\leq& 2e^{-0.0001d(v)}. \nonumber \label{Bvineq} 
\end{eqnarray}
Therefore, by the law of total probability, for every $v\in V$:
\begin{equation}\label{Bvineq}
\mathbf{Pr}\left(B_v\right)  < 2e^{-0.0001d(v)}.
\end{equation}

\subsection{The Local Lemma}

For every $v\in V$ and any event $C_v\in\{A_v,B_v\}$, set 
$$p_{C_v}= p_v: = \frac{1}{1+(d(v))^2}.$$
We define a dependency digraph $D$ with vertex set consisting of all the events $A_v$ and $B_v$, $v\in V$, 
by joining with an arc every $A_v$, and similarly every $B_v$, with all other events $A_u$ and $B_u$ with $u$ at distance at most $2$ from $v$ in $H$. 
(Note that then every event corresponding to a vertex of $D$ is mutually independent of all other considered events which do not belong to its outneighbourhood in $D$.)

Then for every $v\in V$ and $C_v\in\{A_v,B_v\}$, as $\frac{x}{1+x}>e^{-\frac{1}{x}}$ for $x>0$ and $d(u)\geq 2^{-1}d(v)$ for each $u\in N_H(v)$,
\begin{eqnarray}
p_{C_v} \prod_{C'\leftarrow C_v} (1-p_{C'}) 
&\geq& p_v(1-p_v)\prod_{u\in N_H(v)}\prod_{w\in N_H(u)\cup\{u\}\setminus\{v\}}(1-p_w)^2\nonumber\\
&=& \frac{(d(v))^2}{(1+(d(v))^2)^2} \prod_{u\in N_H(v)}
\prod_{w\in N_H(u)\cup\{u\}\setminus\{v\}}\left(\frac{(d(w))^2}{1+(d(w))^2}\right)^2\nonumber\\
&\geq& \frac{(d(v))^2}{(1+(d(v))^2)^2}  \prod_{u\in N_H(v)} \left(\frac{(2^{-1}d(u))^2}{1+(2^{-1}d(u))^2}\right)^{2d(u)} \nonumber\\
&\geq& \frac{(d(v))^2}{(1+(d(v))^2)^2}  \prod_{u\in N_H(v)} e^{-2\cdot 2^2(d(u))^{1-2}} \nonumber\\
&\geq& \frac{(d(v))^2}{(1+(d(v))^2)^2}   e^{-2\cdot 2^2(2^{-1}d(v))^{-1}d(v)} \nonumber\\
&\geq& \frac{1}{2(d(v))^2} e^{-16}. \label{p1-pineq} 
\end{eqnarray}
As for every $v\in V$,
$$d(v) \geq 10^6 \geq 936947 \approx \sqrt[33]{\frac{2\cdot 1920^{70}\cdot e^{16}}{35!}},$$
then by~(\ref{Avineq}) and~(\ref{p1-pineq}), 
\begin{equation}\label{Avineqandp1-pineq}
\mathbf{Pr}(A_v)\leq p_{A_v} \prod_{C'\leftarrow A_v} (1-p_{C'}).
\end{equation}
On the other hand, as $x^2e^{-0.0001x}$ is a decreasing function of $x$ for $x\geq 20000$, while $d(v)\geq 10^6$, we have: 
$$\frac{2e^{-0.0001d(v)}}{ \frac{1}{2(d(v))^2} e^{-16}} = 4e^{16}(d(v))^2e^{-0.0001d(v)}\leq 4\cdot 10^{12}e^{-84}<1,$$
and hence, by~(\ref{Bvineq}) and~(\ref{p1-pineq}),  
\begin{equation}\label{Bvineqandp1-pineq}
\mathbf{Pr}(B_v)\leq p_{B_v} \prod_{C'\leftarrow B_v} (1-p_{C'}).
\end{equation}

Therefore, by~(\ref{Avineqandp1-pineq}), (\ref{Bvineqandp1-pineq}) and the Lov\'asz Local Lemma,
there exist choices of $(c^1_v,c^2_v)\in [0,y_v-1]^2$ for $v\in V$ such that none of the events $A_v, B_v$, $v\in V$ holds, and hence~(\ref{evenG1G2}) holds for every $v\in V$ and $T\leq t$.

\subsection{$\{1,2\}$-weightings of $G_1$ and $G_2$}

We shall now define neighbour sum distinguishing $2$-weightings of $G_1$ and $G_2$. 
First for every vertex $v\in V$ we set 
\begin{equation}\label{lambdavDefinition}
\lambda_v:=4t\cdot y_v, 
\end{equation}
hence by~(\ref{yvBound}),  
\begin{equation}\label{6lambdavIneq}
6\lambda_v \leq 24t \cdot \frac{qd(v)}{24t} = qd(v).
\end{equation}
Next, to every $v\in V$ we assign two sets (corresponding to lists of admissible weighted degrees of $v$ in $G_1,G_2$, resp., modulo $\lambda_v$), where $(j({\rm mod}~2))\in\{0,1\}$ denotes  the remainder of an integer $j$ divided by $2$:
\begin{equation}\label{Av1Definition}
A_v^1:=4tc_v^1+2\left(\log_2 y_v({\rm mod}~2)\right)+\{0,4,8,\ldots,4t-4 \};
\end{equation}
\begin{equation}\label{Av2Definition}
A_v^2:=4tc_v^2+2\left(\log_2 y_v({\rm mod}~2)\right)+\{0,4,8,\ldots,4t-4 \}.
\end{equation}
Let $s_i(v)$ denote the weighted degree (sum) of $v$ with respect to the weighting of $G_i$ we are about to construct; we shall guarantee that the remainder of $s_i(v)$ divided by $\lambda_v$ belongs to $A_v^i\cup(1+A_v^i)$ for $i=1,2$, $v\in V$.
In order to assure sum distinction of (particularly) problematic neighbours, we first specify functions $a_i:V\to \mathbb{Z}$, $i=1,2$, to be utilized within application of Lemma~\ref{1_6Lemma}.
For this aim we start from analyzing all special vertices $v\in V^*$ (cf. ($5^\circ$)) and subsequently for each of these we fix any values
\begin{equation}\label{aiv_assignment1}
a'_i(v)\in A_v^i,
\end{equation} 
$i=1,2$.
Then subsequently for $i=1,2$ and each remaining vertex $v$ ($v\in V\setminus V^*$) one after another, we greedily choose any value 
\begin{equation}\label{aiv_assignment2}
a'_i(v)\in A_v^i
\end{equation} 
which is not yet fixed as $a'_i(u)$ for any neighbour $u$ of $v$ in $H^i_{c_v^1,c_v^2,y_v}$ 
(this is feasible since by our construction, in particular by ($5^\circ$) and the fact $T\leq t$, the degree of $v\in V\setminus V^*$ in $H^i_{c_v^1,c_v^2,y_v}$ is at most $t-1$, while by~(\ref{Av1Definition}) and~(\ref{Av2Definition}), $|A_v^i|=t$). Finally, we set: 
\begin{equation}\label{a'iv_assignment}
a_i(v):=a'_i(v)-d_{G_i}(v)
\end{equation} 
for each $i=1,2$, $v\in V$.
Then for $i=1,2$ we apply Lemma~\ref{1_6Lemma} to the graph $G_i$ with $a_i:V\to \mathbb{Z}$ defined in~(\ref{a'iv_assignment}) and $\lambda_v$ defined in~(\ref{lambdavDefinition}), what is feasible due to~(\ref{evenG1G2}) and~(\ref{6lambdavIneq}), in order to obtain a subgraph $H'''_i$ of $G_i$ with 
\begin{equation}\label{dH''iIn}
d_{H'''_i}(v)\in\left[\frac{d_{G_i}(v)}{3},\frac{2d_{G_i}(v)}{3}\right],
\end{equation}
\begin{equation}\label{dH''iEquiv}
d_{H'''_i}(v)\equiv {a_i(v),a_i(v)+1} \pmod {\lambda_v}
\end{equation}
for each $v\in V$, and afterwards we define a $\{1,2\}$-weighting $c_i$ of $G_i$ by setting:
$$c_i(e)=\left\{\begin{array}{lll} 
2 &{\rm if} & e\in E(H'''_i),\\
1 &{\rm if} & e\in E(G_i)\setminus E(H'''_i).
\end{array}\right.$$
Note that then, by~(\ref{a'iv_assignment}) and~(\ref{dH''iEquiv}),
\begin{eqnarray}
s_i(v) &=& 2d_{H'''_i}(v)+1(d_{G_i}(v)-d_{H'''_i}(v)) = d_{G_i}(v)+d_{H'''_i}(v) \label{sidGdH}\\
&\equiv& {a'_i(v),a'_i(v)+1} \pmod {\lambda_v} \label{sidGdHmodLambdav}
\end{eqnarray}
for $i=1,2$, $v\in V$.

\subsection{Final Analysis}

Fix $i\in\{1,2\}$ and suppose $e=uv\in E(G_i)$. Without loss of generality we assume that $$d(u)\geq d(v).$$

If $e\notin H_i$, hence $d(u)> 2d(v)$,
then by~(\ref{evenG1G2}), (\ref{dH''iIn}) and~(\ref{sidGdH}),
\begin{eqnarray}s_i(u)&=& d_{G_i}(u)+d_{H'''_i}(u) \geq \frac{4}{3}d_{G_i}(u) \geq \frac{4}{3} q d(u)
> \frac{8}{3} q d(v) >  \frac{5}{3} (1-q) d(v)\nonumber\\
&\geq&  \frac{5}{3} \left(d(v)-d_{G_{3-i}}(v)\right) = \frac{5}{3} d_{G_i}(v) \geq d_{G_i}(v)+d_{H'''_i}(v) = s_i(v).\nonumber
\end{eqnarray}
Note on the other hand that if $e\in H_i$, then by~(\ref{yuyvBound}), $y_u=y_v$ or $y_u=2y_v$.

If $e\in H_i$ and $y_u=2y_v$, then by~(\ref{Av1Definition})--(\ref{aiv_assignment2}) and~(\ref{sidGdHmodLambdav}), 
$$s_i(v)\equiv 2\log_2y_v,2\log_2y_v+1 \pmod{4}~~~~~{\rm and}~~~~~s_i(u)\equiv 2\log_2y_v+2,2\log_2y_v+3 \pmod{4}.$$

If $e\in H_i$, $y_u=y_v$ and $(c^1_v,c^2_v)\neq (c^1_u,c^2_u)$, then by~($1^\circ$) and~($2^\circ$), $c^i_u\neq c^i_v$, and hence $A_u^i$ and $A_v^i$ are disjoint sets (of even integers), and hence by~(\ref{Av1Definition})--(\ref{aiv_assignment2}) and~(\ref{sidGdHmodLambdav}),
$$ s_i(u) \not\equiv s_i(v) \pmod {\lambda_v}.$$

If finally $e\in H_i$, $y_u=y_v$ and $(c^1_v,c^2_v) = (c^1_u,c^2_u)$, then $e\in E(H^i_{c_v^1,c_v^2,y_v})$, and hence by our construction of $a'_i$, we conclude that $a'_i(u)$ and $a'_i(v)$ are two distinct even integers in $[0,\lambda_v-2]$, and thus by~(\ref{sidGdHmodLambdav}),
$$ s_i(u) \not\equiv s_i(v) \pmod {\lambda_v}.$$

In all cases we thus obtain that $s_i(u)\neq s_i(v)$, and hence $G_i$ is $\{1,2\}$-weight colourable for $i=1,2$. \qed

\section{Remarks}

Though we have put some effort into designing our proof and appropriate selection of constants used, 
our approach could still be optimized in several aspects. This might however influence the clarity of presentation of the argument, and it is doubtful that applying the general proving scheme we propose one could show Theorem~\ref{MainTh22Standard} to be valid e.g. for $\delta\geq 10^5$. This constant could be further improved in the case of regular graphs though, by means of the symmetric version of the Lov\'asz Local Lemma, but presumably not below  
$10^4$, we thus omit details.

\end{document}